\documentclass{article}

\usepackage[english]{babel}
\usepackage{latexsym}
\usepackage{graphicx}
\usepackage{amsmath}
\usepackage{epsfig}
\usepackage{amssymb}
\usepackage{mathrsfs}
\usepackage{amsthm}

\def\R{{\mathbb{R}}}

\def\1{\mbox{$\mathbf{1}$}}

\def\det{\mathop {\rm det}\nolimits}

\def\div{\mathop{\rm div}\nolimits}

\def\diam{\mathop{\rm diam}\nolimits}

\newtheorem{defi}{\indent Definition}[section]
\newtheorem{teo}[defi]{\indent Theorem}
\newtheorem{lemma}[defi]{\indent Lemma}

\newtheorem{prop}[defi]{\indent Proposition}

\newtheorem{oss}[defi]{\indent Remark}

\numberwithin{equation}{section}

\title{\ Determining the anisotropic traction state in a membrane \\
by boundary
measurements}

\author{{ G. Alessandrini  }
\thanks{Dipartimento di Matematica e Informatica,
Via Valerio, 12 - 34127 Trieste, Italy, email: {\tt alessang@univ.trieste.it}. Work supported in part by MIUR, PRIN n. 2004011204 and by GNAMPA,
INdAM Progetto \emph{Problemi al contorno inversi} 2006 .}\and { E. Cabib}
\thanks{Dipartimento di Ingegneria Civile,
Via delle Scienze, 208 - 33100 Udine, Italy, email: {\tt
cabib@uniud.it} }}

\date{}

\begin{document}

\maketitle


\section{Introduction}

\hspace{\parindent} Consider an elastic thin membrane which
occupies a planar region represented by a simply connected bounded
open set $\Omega\subset\R^2$. In the same plane a vector force
field $T$ is applied on its boundary $\partial\Omega$, so that the
membrane is subject to a distributed {\it pretraction state},
expressed by a positive and symmetric tensor
$\sigma=\{\sigma_{ij}\}$, $i,j=1,2$, which satisfies the plane
equilibrium equations
\begin{eqnarray}\label{traction eq}
                 &\div\sigma=0 \ , \  \mbox{in} \ \Omega \ ,
                 \hfill&\\
                \label{traction eq bordo} & \sigma \nu=T\ ,\ \mbox{on}\ \partial\Omega \ ,
\end{eqnarray}
where $\nu$ is the outer unit normal to $\partial\Omega$. The
transverse displacement $u(x)$ of the membrane will be governed by
the equation
\begin{eqnarray}
                & \label{ellipeq}-\div(\sigma \nabla u)=f \ , \ \mbox{in} \ \Omega \ ,&
                 \\
                & \label{ellipeqbordo}\quad u=\varphi \ ,\ \mbox{on}\ \partial\Omega \ ,&
\end{eqnarray}
where $f(x)$  represents the distributed transverse load applied
to it, and $\varphi$ represents the prescribed transverse
displacement at the boundary.

In this note we wish to investigate the inverse problem of
determining the plane traction state tensor $\sigma$ from boundary
measurements on displacements and on the corresponding forces. As
an initial attempt, since the external load $f$ has no influence
on the pretraction state $\sigma$, and in analogy with many other
well-known inverse boundary problems, see for instance
\cite{a1988,c1980,nu1994}, it seems natural to treat the case when
$f=0$ in \eqref{ellipeq}, and consider as available data an
arbitrary transverse displacement $\varphi$ on the boundary and
the corresponding transverse load on the boundary, namely the
reaction of the boundary constraints
\begin{equation}\label{dirichneum}
\Lambda_{\sigma} \varphi = \sigma \nabla u \cdot \nu \ .
\end{equation}

Hence, fixing any $\varphi \in H^{1/2}(\partial \Omega)$, if we denote by $u\in H^1(\Omega)$ the weak solution to the Dirichlet problem

\begin{eqnarray}
                & \label{equazomog}\div(\sigma \nabla u)=0 \ ,  \mbox{in} \ \Omega \ ,&
                 \\
                & \label{datoDirichlet}\quad u=\varphi \ , \mbox{on}\ \partial\Omega \ ,&
\end{eqnarray}
we introduce the Dirichlet-to-Neumann (D-N) map as the bounded linear operator
\begin{equation}\label{dirichneummap}
\Lambda_{\sigma} : H^{\frac{1}{2}}(\partial \Omega) \rightarrow H^{-\frac{1}{2}}(\partial \Omega)\ ,
\end{equation}
defined, in weak terms, by the formula
\begin{equation}\label{dirichneumweak}
<\Lambda_{\sigma} \varphi, v_{|{\partial \Omega}}> =
\int_{\Omega}\sigma \nabla u \cdot \nabla v \ ,\quad \text{ for
every } v \in H^1(\Omega)\ ,
\end{equation}
where $v_{|\partial \Omega}$ denotes the trace on $\partial \Omega$ of $ v \in H^1(\Omega) $.

Thus we examine here the problem of determining $\sigma$ from the
knowledge of $\Lambda_{\sigma}$. The peculiarity of this problem
is that, by its own nature, the tensor $\sigma$ is anisotropic,
and since Tartar's example,  as reported in \cite{kv1984}, it is
well-known that a general anisotropic tensor $\sigma$ cannot be
uniquely determined by the D-N map $\Lambda_{\sigma}$. In this
case however, we shall see, in the next Section \ref{unique} that
we can take advantage of the null divergence condition
\eqref{traction eq} and obtain the uniqueness in such restricted
class of genuinely anisotropic tensors, see Theorem \ref{teorema}.

In Section \ref{stab} we also consider the question of stability,
that is of the continuous dependence upon the data. And although
the present results are very preliminary, they show up interesting
phenomena, which are markedly different from those available  for
the well-known \emph{inverse conductivity problem}. In fact, on
one hand, we prove, Theorem \ref{stability}, a qualitative form of
stability, when a very weak topology is assigned on the class of
tensors, namely the topology of \emph{G-convegence}. And, on the
other hand, we show that a Lipschitz stability bound holds for the
mean value of $\sigma$.

\section{Uniqueness}\label{unique}

\hspace{\parindent}For any given $K\geq 1$ we consider the class of
tensors
\begin{equation}\label{MK} M_K=\{\sigma\in L^\infty(\Omega,M^{2\times
2})\;|\;K^{-1}|\xi|^2\leq \sigma\xi\cdot\xi\leq K |\xi|^2\quad \text{for every } \xi\in\R^2 \ \}\ ,\end{equation}
here $M^{2\times 2}$ denotes
the set of ${2\times 2}$ symmetric matrices. Let us also introduce
\begin{equation}\label{SigmaK}\Sigma_K=\{\sigma\in M_K\;|\;\div\sigma=0\}\ ,\end{equation} where the null divergence condition $\div\sigma=0$ is meant in
the weak sense
\begin{equation}\label{weakdivfree}
\int_\Omega\sigma \nabla v =0 \ \quad\text{for every }   v \in
H_0^1(\Omega)\ .
\end{equation}
Let us also introduce
\begin{equation}\label{Sigma} M = \cup_{K\geq 1}M_K \ , \ \Sigma = \cup_{K\geq 1}\Sigma_K \  . \end{equation}
Given a $W^{1,2}$ mapping $\Phi : {\Omega} \rightarrow {D} \subset \mathbb{R}^2$ we denote its Jacobian matrix as follows
\begin{equation}\label{jac}
D \Phi (x) = \left\{\frac{\partial \Phi_i(x)}{\partial{x_j}}\right\}\hspace{.5truecm}i,j=1,2 \ .
\end{equation}
We recall that $\Phi$ is said to be \emph{quasiconformal} if, for some $Q\geq 1$ it satisfies
\begin{equation}\label{QC}
\|D \Phi\|^2 \leq Q \det D \Phi \text{  a.e. in } \Omega\ ,
\end{equation}
and it is invertible, see for instance Ahlfors \cite{ah1966}.
Here, for any matrix $A$, we denote $\|A\|^2=\text{tr}AA^T$ and
the suffix $T$ denotes transpose.

For any tensor $\sigma \in M_K$ and any quasiconformal mapping $\Phi$, we introduce
\begin{equation}\label{pushf}
T_{\Phi}\sigma(y) = \frac{D\Phi \sigma D \Phi^T}{\det D \Phi}(\Phi^{-1}(y)) \ ,\quad \text{ for every } y \in D \ .
\end{equation}
 This new tensor defined in
$D$, called the {\it push-forward} of $\sigma$ by $\Phi$, is again symmetric and satisfies the ellipticity condition with some possibly new
$K\geq 1$. Moreover, one can verify that such operation preserves the bilinear Dirichlet form associated to $\sigma$, that is
\begin{equation}\label{forminvariance}
\int_{\Omega}\sigma \nabla u \cdot \nabla v = \int_D T_{\Phi}\sigma
\nabla (u \circ \Phi^{-1}) \cdot \nabla (v \circ \Phi^{-1}) \
\quad\text{ for every } u,v \in H^1(\Omega) \ .
\end{equation}

\begin{teo}\label{teorema}$\Lambda_\sigma$ uniquely determines $\sigma$ among all
tensors in $\Sigma$.
\end{teo}

The proof will be a consequence of the following two Lemmas.
\begin{lemma}\label{lemmatwisttensor} Let $\sigma\in \Sigma$. Then the tensor
\begin{equation}\label{twisttensor}
T_{\Phi}\sigma (D \Phi^{-1})^T
\end{equation}
is divergence free in $D$.
\end{lemma}
\begin{proof}
Condition \eqref{weakdivfree} can be rewritten as
\begin{equation}\label{weakdivfreeindex}
\int_\Omega\sigma \nabla x_j \cdot \nabla v =0 \ \quad\text{ for
every } v \in H_0^1(\Omega)\  , \ j=1,2 \ .
\end{equation}
By \eqref{forminvariance} one computes
\begin{equation}\label{weakdivfreeindexy}
\int_D T_{\Phi}\sigma \nabla ((\Phi^{-1})_j) \cdot \nabla (v \circ
\Phi^{-1}) =0 \ \text{ for every }  v \in H_0^1(\Omega)\ , \ j=1,2 \
,
\end{equation}
and the thesis follows.
\end{proof}
\begin{lemma}\label{lemmarevtwisttensor}
Let $\sigma\in M$ and suppose that for a given mapping  $\Phi$ we have $\div T_{\Phi}\sigma=0$ in $D$. Then
\begin{equation}\label{revdivfree}
\div(\sigma D \Phi^T)=0 \ .
\end{equation}
\end{lemma}
\begin{proof}
The proof follows immediately from Lemma \ref{lemmatwisttensor},
just by reversing the roles of $\Phi$ , $\Phi^{-1}$ and of
$\sigma$, $T_{\Phi}\sigma$, respectively. \end{proof}

We are now in a position to prove our main result.
\begin{proof}[Proof of Theorem \ref{teorema}]
By the results of Astala,   P\"{a}iv\"{a}rinta and Lassas
\cite[Theorem 1]{alp2005}, which have extended to the $L^{\infty}$
setting those of Sylvester \cite{s1990} and Nachman \cite{n1996},
we have that $\Lambda_\sigma$ determines uniquely the class
\begin{equation}\label{class}E_{\sigma}=\{\sigma'\in M | \sigma'=T_{\Phi}\sigma, \text{ with } \Phi:\Omega\rightarrow\Omega \text{ quasiconformal and such
that }\Phi_{|\partial\Omega}=I\}  \ . \end{equation}
The class $E_{\sigma}$ contains at most one divergence free element. In fact, if $\Phi$ is
a quasiconformal mappings which fixes the boundary, and such that $ \div T_{\Phi}\sigma = 0 $ then, by Lemma \ref{lemmarevtwisttensor}, we have
\begin{equation}\label{eq: identity}\left\{\begin{array}{*{20}c}
                           \div(\sigma D \Phi^T)=0\hfill & \mbox{in}\;\Omega \ ,\hfill\\
                           \Phi=I\hfill &
                           \mbox{on}\;\partial\Omega\ .\hfill
                           \end{array}\right.\end{equation}
Note that this system, is formed by two uncoupled  Dirichlet problems for the two components of the mapping $\Phi$. On the other hand we observe
that, if $\sigma$ is divergence free, then the identity mapping $I$ is itself a solution to (\ref{eq: identity}) and by uniqueness for the
Dirichlet problem, we obtain $\Phi=I$ on $\Omega$.
\end{proof}

\begin{oss}\label{rem1}
It is worth mentioning, that the crucial fact used in this proof
is that all linear functions are solutions to the elliptic
equation \eqref{ellipeq}, this is a condition on $\sigma$ which is
in fact equivalent to the null divergence condition
\eqref{weakdivfree}. Indeed this property has been used already in
a study on optimization of tension structures in the different
context of variational problems and G-convergence, see
\cite{cdr1990}, \cite{cd1991} and also Section \ref{stab} below.

A further application of this property of linear functions is the
possibility to identify the traction $T$ applied on the boundary.
For every $\xi\in\R^2$ let $\varphi_\xi(x)=\xi\cdot x$ be the
Dirichlet data. Since the corresponding solution is
$u_\xi(x)=\xi\cdot x$ all over $\Omega$, we have
\begin{equation}\label{bdrytraction}
\Lambda_\sigma\varphi_\xi=\sigma\nabla
u_\xi\cdot\nu=\sigma\xi\cdot\nu=\sigma\nu\cdot\xi=T\cdot\xi\,,
\end{equation}
whose knowledge for every $\xi\in\R^n$ is equivalent to the
knowledge of $T$.

If the same argument is applied to the particular case of a square
network $Q=[0,1]\times[0,1]$, that is, a portion of fabrics made
by two families of parallel elastic strings which cross
orthogonally, the identification of $\sigma$ is immediate. The
particular situation leads to define as admissible all the tensors
of the form
$$\sigma(x)=\left(\begin{array}{*{20}c}
           \sigma_1(x_2) & 0\\
           0 & \sigma_2(x_1)\end{array}\right)\,.$$
with $K^{-1}\leq \sigma_1,\sigma_2\leq K$. Consider the Dirichlet
data $\varphi(x)=x_1$ on $\partial Q$, so that also $u(x)=x_1$ on
$Q$, and the corresponding $\psi(x_2)=\Lambda_\sigma\varphi$ on
the edge $x_1=1$ where $\nu=(1,0)$. Then we have
$$\sigma_1(x_2)=\sigma(1,x_2)\nabla u\cdot\nu=\psi(x_2)\, ,$$
likewise, $\sigma_2(x_1)$ can be identified as well.
\end{oss}

\section{Stability}\label{stab}

\hspace{\parindent}This section is devoted to the continuity
properties of the inverse of the map $$ \Sigma_K \ni
\sigma\rightarrow\Lambda_\sigma \in \mathscr
L(H^{1/2}(\Omega),H^{-1/2}(\Omega))$$ when we assign to $\Sigma_K$
the topology of  G-convergence. Let us recall here the basic
notions  and some important properties of the G-convergence. A
wide literature is available on this subject, we refer for example
to the classical papers \cite{ds1973,ms1969,s1968,s1975} and to
the book by \mbox{Dal Maso} \cite{d1993} where G-convegence is
cast in the more general theory of $\Gamma$-convegence.

\begin{defi}\label{g-conv} A sequence $\{\sigma_h\}\subset M_K$ is said to G-converge to
$\sigma\in M_K$, and we write $\sigma_h\stackrel{G}{\to} \sigma$,
if for every $f\in H^{-1}(\Omega)$ the corresponding sequence
$\{u_h\}\subset H^1_0(\Omega)$ of solutions to the
\emph{inhomogeneous} problems
\begin{equation}\label{g-succ}
-\div(\sigma_h\nabla u_h)=f\quad \ \textrm{in} \ \Omega \ , \
u_h=0\quad \  \textrm{on} \ \partial \Omega \ ,
\end{equation}
converges weakly in $H^1_0(\Omega)$ to the solution $u\in
H^1_0(\Omega)$ of the problem
\begin{equation}\label{g-lim}
-\div(\sigma\nabla u)=f\quad  \ \textrm{in} \ \Omega\,\ , \
u=0\quad \textrm{on} \
\partial \Omega \ .
\end{equation}
\end{defi}
It is well known that G-convergence is induced by a compact
metrizable topology on $M_K$, \cite[Remark 4]{s1975}.

It is also worth recalling that the $L^1_{\rm loc}$-strong
convergence implies the G-convergence, \cite[Proposition
5]{s1968}, \cite[Remark 11]{s1975}.

\begin{teo}\label{stability} Given $K\geq1$, the mapping $ \Sigma_K \ni \sigma \rightarrow\Lambda_\sigma \in \mathscr L(H^{1/2}(\Omega),H^{-1/2}(\Omega))$ has a
continuous inverse when $\Sigma_K$ is endowed with the topology of
G-convergence.
\end{teo}
\begin{proof}
First we recall that \cite{cd1991}, by the characterization of
$\Sigma_K$ as the subclass of those $\sigma \in M_K$ for which all
linear functions are solutions to \eqref{equazomog}, implies that
$\Sigma_K$ is a closed set in the G-topology and hence it is
compact. Let $\{\sigma_h\}\subset\Sigma_K$ and $\sigma \in
\Sigma_K$ be such that
$\|\Lambda_{\sigma_h}-\Lambda_\sigma\|\rightarrow 0$. By the above
mentioned  compactness , there exists a subsequence
$\{\sigma_{r_h}\}$ of $\{\sigma_h\}$ such that $\sigma_{r_h}
\stackrel{G}{\rightarrow}\sigma'\in\Sigma_K$ and we prove
$\sigma'=\sigma$. For any $\varphi\in H^{1/2}(\partial \Omega)$,
let $u_h$, $u'$ the solutions to \eqref{equazomog} when $\sigma$
is replaced with with $\sigma_h$, $\sigma'$, respectively. Then,
by convergence of the energies \cite{s1975}, we have
\begin{equation}<\Lambda_{\sigma_{r_h}}\varphi,\varphi>= \int_{\Omega}\sigma_{r_h}\nabla u_{r_h}\cdot \nabla
u_{r_h}\rightarrow \int_{\Omega}\sigma'\nabla u'\cdot \nabla u'=
<\Lambda_{\sigma'}\varphi,\varphi> \ ,
\end{equation}
on the other hand, $\Lambda_{\sigma_h}\rightarrow\Lambda_{\sigma}$,
therefore $<\Lambda_{\sigma'}\varphi,\varphi>=
<\Lambda_{\sigma}\varphi,\varphi>$ for every $\varphi\in
H^{1/2}(\partial \Omega)$. From the uniqueness Theorem
\ref{teorema}, we get $\sigma'=\sigma$. The above argument applies
to any subsequence of $\{\sigma_h\}$. Thus we have obtained that for
any subsequence of $\{\sigma_h\}$ there is a sub-subsequence which
G-converges to $\sigma$, and hence, the full sequence $\{\sigma_h\}$
must G-converge to $\sigma$.
\end{proof}

\begin{oss}\label{confr-iso} It is worthwhile to compare this
result with the case of the inverse conductivity problem, that is
when the unknown $\sigma\in M_K$ is a-priori known to be
isotropic, that is $\sigma= \gamma I$ where $\gamma\in
L^{\infty}(\Omega)$ is a scalar function satisfying
$K^{-1}\leq\gamma\leq K$ and $I$ denotes the identity matrix.
Indeed for the inverse conductivity problem, stability with
respect to G-convergence fails, see \cite{kv1987} for  related
arguments. In fact, as is well-known, Marino and Spagnolo
\cite{ms1969} proved that there exist a constant $c>1$, depending
only on the space dimension $n$ (in our case $n=2$), such that any
tensor in $M_{K/c}$ can be approximated in the sense of
G-convegence by isotropic tensors of $M_K$. Hence if we had
stability with respect to G-convergence for isotropic tensors,
that would imply the uniqueness in the class $M_{K/c}$ of
\emph{anisotropic} tensors, which, by the above mentioned example
of Tartar cannot hold true. Hence the stability result above,
Theorem~\ref{stability},  is crucially based on the property of
our set of admissible matrices, $\Sigma_K$, of being G-closed.
\end{oss}

We recall that in  \cite{cdr1990}, and in \cite{cd1991} in a more
general context, it was  proved that on $\Sigma_K$  the
G-convergence is equivalent to the $L^{\infty}(\Omega)$-weak*
convergence. Therefore, as a consequence of Theorem
\ref{stability}, we also obtain that for every $\psi \in
L^{1}(\Omega)$ and for every $i,j=1,2$ the functional $F$ defined
by
\begin{equation}\label{functional}
F(\sigma)=\int_{\Omega}\psi \sigma_{i,j} \ , \ \sigma\in \Sigma_K
\ ,
\end{equation}
depends continuously on $\Lambda_{\sigma}$.

In the very special case when, in \eqref{functional}, we choose
$\psi \equiv \frac{1}{|\Omega|}$ a concrete stability estimate can
be obtained. In fact, in the next Proposition we  show that the
average of $\sigma$, a quantity which can be interpreted as a
global measure of the pretraction field, depends in a Lipschitz
continuous fashion  on the Dirichlet-to-Neumann map.
\begin{prop} For any  $\sigma, \sigma' \in \Sigma$ we have
\begin{equation} \left\|\frac{1}{|\Omega|}\int_{\Omega}(\sigma-\sigma')\right\| \leq (1+(\diam\Omega)^2)\|\Lambda_\sigma-\Lambda_{\sigma'}\| \ .
\end{equation}
\end{prop}
\begin{proof} Being linear functions solutions, we can use \eqref{dirichneumweak}, with $u=x_i \ , v=x_j \ , i,j=1,2$, both for $\sigma$ and
$\sigma'$. We obtain
\begin{equation} \int_{\Omega}(\sigma-\sigma')_{ij} = <(\Lambda_\sigma-\Lambda_{\sigma'})x_i,x_j> \ ,
\end{equation}
and the thesis follows by straightforward computations.\end{proof}
\begin{oss}\label{aver}
Also in this case it may be interesting to make a comparison with
the inverse conductivity problem. In fact it is an open problem
whether, for the average $\frac{1}{|\Omega|}\int_{\Omega}\gamma$
of an isotropic tensor $\sigma= \gamma I$, the  Lipschitz
stability in terms of the corresponding the Dirichlet-to-Neumann
map holds true, see \cite{av2005}.
\end{oss}

\renewcommand*{\refname}{References}

\end{document}